\documentclass{amsart}
\usepackage{amssymb}

\usepackage{amscd}
\usepackage{thmdefs}

\input{tcilatex}
\theoremstyle{definition}
\theoremstyle{remark}
\numberwithin{equation}{section}

\input tcilatex

\begin{document}
\title[Crossed Product Probability Spaces]{Group-Freeness and Certain Amalgamated Freeness}
\author{Ilwoo Cho}
\address{Saint Ambrose Univ., Dep of Math, 116 McM hall, 158 W. Locust St.,
Davenport, Iowa 52803, U. S. A.}
\email{chowoo@sau.edu}
\thanks{The author specially thanks to Prof. F. Radulescu, who is his Ph. D. thesis
advisor in Univ. of Iowa, for the valuable discussion and advice. Also, the
author appreciate all supports from St. Ambrose Univ.. In particular, the
author thanks to Prof. V. Vega and Prof. T. Anderson, for the useful
discussion and for the kind encouragement and advice.}
\date{03 / 28 / 2006}
\dedicatory{}
\thanks{}
\maketitle

\begin{abstract}
In this paper, we will consider certain amalgamated free product structure
in crossed product algebras. Let $M$ be a von Neumann algebra acting on a
Hilbert space $H$ and $G,$ a group and let $\alpha $ $:$ $G$ $\rightarrow $ $%
AutM$ be an action of $G$ on $M$, where $AutM$ is the group of all
automorphisms on $M.$ Then the crossed product $\Bbb{M}$ $=$ $M$ $\times
_{\alpha }$ $G$ of $M$ and $G$ with respect to $\alpha $ is a von Neumann
algebra acting on $H$ $\otimes $ $l^{2}(G),$ generated by $M$ and $%
\{u_{g}\}_{g\in G},$ where $u_{g}$ is the unitary representation of $g$ on $%
l^{2}(G).$ We show that $M$ $\times _{\alpha }(G_{1}$ $*$ $G_{2})$ $=$ $(M$ $%
\times _{\alpha }$ $G_{1})$ $*_{M}$ $(M$ $\times _{\alpha }$ $G_{2}).$ We
compute moments and cumulants of operators in $\Bbb{M}.$ By doing that, we
can verify that there is a close relation between Group Freeness and
Amalgamated Freeness under the crossed product. As an application, we can
show that if $F_{N}$ is the free group with $N$-generators, then the crossed
product algebra $L_{M}(F_{n})$ $\equiv $ $M\ \times _{\alpha }$ $F_{n}$
satisfies that $L_{M}(F_{n})$ $=$ $L_{M}(F_{k_{1}})$ $*_{M}$ $%
L_{M}(F_{k_{2}}),$ whenever $n$ $=$ $k_{1}$ $+$ $k_{2},$ for $n,$ $k_{1},$ $%
k_{2}$ $\in $ $\Bbb{N}.$
\end{abstract}

\strut

In this paper, we will consider a relation between a free product of groups
and a certain free product of von Neumann algebras with amalgamation over a
fixed von Neumann subalgebra. In particular, we observe such relation when
we have crossed product algebras. Crossed product algebras have been studied
by various mathematicians. Let $M$ be a von Neumann algebra acting on a
Hilbert space $H$ and $G,$ a group, and let $\Bbb{M}$ $=$ $M$ $\times
_{\alpha }$ $G$ be the crossed product of $M$ and $G$ via an action $\alpha $
$:$ $G$ $\rightarrow $ $AutM$ of $G$ on $M,$ where $AutM$ is the
automorphism group of $M.$ This new von Neumann algebra $\Bbb{M}$ acts on
the Hilbert space $H$ $\otimes $ $l^{2}(G),$ where $l^{2}(G)$ is the group
Hilbert space. Each element $x$ in $\Bbb{M}$ has its Fourier expansion

\strut

\begin{center}
$x=\underset{g\in G}{\sum }m_{g}u_{g},$ \ for \ $m_{g}$ $\in $ $M$
\end{center}

\strut

where $u_{g}$ is the (left regular) unitary representation of $g$ $\in $ $G$
on $l^{2}(G).$

\strut

On $\Bbb{M},$ we have the following basic computations;

\strut

(0.1) If $u_{h}$ is the unitary representation of $h$ $\in $ $G,$ as an
element in $\Bbb{M},$ then

$\strut $

\begin{center}
$u_{g_{1}}u_{g_{2}}=u_{g_{1}g_{2}}$ \ \ and \ \ $u_{g}^{*}=u_{g^{-1}},$ for
all $g,g_{1},g_{2}\in G$
\end{center}

\strut

(0.2) \ If $m_{1},$ $m_{2}$ $\in $ $M$ and $g_{1},$ $g_{2}$ $\in $ $G,$ then

$\strut $

\begin{center}
$
\begin{array}{ll}
\left( m_{1}u_{g_{1}}\right) \left( m_{2}u_{g_{2}}\right) & 
=m_{1}u_{g_{1}}m_{2}(u_{g_{1}}^{-1}u_{g_{1}})u_{g_{2}} \\ 
&  \\ 
& =\left( m_{1}\left( \alpha _{g_{1}}(m_{2})\right) \right) u_{g_{1}g_{2}}
\end{array}
$
\end{center}

\strut

(0.3) If $mu_{g}\in \Bbb{M},$ then

\strut

\begin{center}
$
\begin{array}{ll}
\left( mu_{g}\right) ^{*} & =u_{g}^{*}m^{*}=u_{g^{-1}}m^{*}\left(
u_{g}u_{g^{-1}}\right) \\ 
&  \\ 
& =\left( \alpha _{g^{-1}}(m^{*})\right) u_{g^{-1}}=\left( \alpha
_{g^{-1}}(m^{*})u_{g}^{*}\right) .
\end{array}
$
\end{center}

\strut \strut

(0.4) If $m\in M$ and $g\in G,$ then

\strut

\begin{center}
$u_{g}m=u_{g}mu_{g^{-1}}u_{g}=\alpha _{g}(m)u_{g}.$
\end{center}

\strut

and

\begin{center}
$mu_{g}=u_{g}u_{g^{-1}}mu_{g}=u_{g}\cdot \alpha _{g^{-1}}(m)$
\end{center}

\strut \strut

The element $u_{g}m$ is of course contained in $\Bbb{M}$, since it can be
regarded as $u_{g}mu_{e_{G}},$ where $e_{G}$ is the group identity of $G,$
for $m$ $\in $ $M$ and $g$ $\in $ $G.$

\strut

Free Probability has been researched from mid 1980's. There are two
approaches to study it; the Voiculescu's original analytic approach and the
Speicher's combinatorial approach. We will use the Speicher's approach. Let $%
M$ be a von Neumann algebra and $N,$ a $W^{*}$-subalgebra and assume that
there is a conditional expectation $E$ $:$ $M$ $\rightarrow $ $N$ satisfying
that (i) $E$ is a continous $\Bbb{C}$-linear map, (ii) $E(n)$ $=$ $n,$ for
all $n$ $\in $ $N$, (iii) $E(n_{1}$ $m$ $n_{2})$ $=$ $n_{1}$ $E(m)$ $n_{2},$
for all $m$ $\in $ $M$ and $n_{1},$ $n_{2}$ $\in $ $N,$ and (iv) $E(m^{*})$ $%
=$ $E(m)^{*},$ for all $m$ $\in $ $M.$ If $N$ $=$ $\Bbb{C},$ then $E$ is a
continous linear functional on $M,$ satisfying that $E(m^{*})$ $=$ $%
\overline{E(m)},$ for all $m$ $\in $ $M.$ The algebraic pair $(M,$ $E)$ is
called an $N$-valued $W^{*}$-probability space. All operators $m$ in $(M,$ $%
E)$ are said to be $N$-valued random variables. Let $x_{1},$ ..., $x_{s}$ $%
\in $ $(M,$ $E)$ be $N$-valued random variables, for $s$ $\in $ $\Bbb{N}.$
Then $x_{1},$ ..., $x_{s}$ contain the following free distributional data.

\strut

$\circ $ \ $(i_{1},$ ..., $i_{n})$-th joint $*$-moment : \ $E\left(
x_{i_{1}}^{u_{i_{1}}}...x_{i_{n}}^{u_{i_{n}}}\right) $

\strut

$\circ $ \ $(j_{1},$ ..., $j_{m})$-th joint $*$-cumulant : $k_{m}\left(
x_{j_{1}}^{u_{j_{1}}},\text{ ..., }x_{j_{m}}^{u_{j_{m}}}\right) $ such that

\strut

\begin{center}
$k_{m}\left( x_{j_{1}}^{u_{j_{1}}},\text{ ..., }x_{j_{m}}^{u_{j_{m}}}\right) 
\overset{def}{=}\underset{\pi \in NC(m)}{\sum }E_{\pi }\left(
x_{j_{1}}^{u_{j_{1}}},\text{ ..., }x_{j_{m}}^{u_{j_{m}}}\right) \mu (\pi
,1_{m}),$
\end{center}

\strut

for $(i_{1},$ ..., $i_{n})$ $\in $ $\{1,$ ..., $s\}^{n}$, $(j_{1},$ ..., $%
j_{m})$ $\in $ $\{1,$ ..., $s\}^{m},$ for $n,$ $m$ $\in $ $\Bbb{N},$ and $%
u_{i_{k}},$ $u_{j_{i}}$ $\in $ $\{1,$ $*\},$ and where $NC(m)$ is the
lattice of all noncrossing partitions over $\{1,$ ..., $m\}$ with its
minimal element $0_{m}$ $=$ $\{(1),$ ..., $(m)\}$ and its maximal element $%
1_{n}$ $=$ $\{(1,$ ..., $m)\}$ and $\mu $ is the M\"{o}bius functional in
the incidence algebra and $E_{\pi }(...)$ is the partition-depending moment
of $x_{j_{1}},$ ..., $x_{j_{m}}$ (See [19]).

\strut

For instance, $\pi $ $=$ $\{(1,$ $4),$ $(2,$ $3)\}$ is in $NC(4).$ We say
that the elements $(1,$ $4)$ and $(2,$ $3)$ of $\pi $ are blocks of $\pi ,$
and write $(1,$ $4)$ $\in $ $\pi $ and $(2,$ $3)$ $\in $ $\pi .$ In this
case, the partition-depending moment $E_{\pi }(x_{j_{1}},$ ..., $x_{j_{4}})$
is determined by

\strut

\begin{center}
$E_{\pi }\left( x_{j_{1}},x_{j_{2}},x_{j_{3}},x_{j_{4}}\right) =E\left(
x_{j_{1}}E(x_{j_{2}}x_{j_{3}})x_{j_{4}}\right) .$
\end{center}

\strut

The ordering on $NC(m)$ is defined by

\strut

\begin{center}
$\pi \leq \theta \Longleftrightarrow $ for any block $B$ $\in $ $\pi ,$
there is $V$ $\in \theta $ such that $B$ $\subseteq $ $V,$
\end{center}

\strut

for $\pi ,$ $\theta $ $\in $ $NC(m),$ where ``$\subseteq $'' means the usual
set-inclusion.

\strut \strut

Suppose $M_{1}$ and $M_{2}$ are $W^{*}$-subalgebras of $M$ containing their
common subalgebra $N.$ The $W^{*}$-subalgebras $M_{1}$ and $M_{2}$ are said
to be free over $N$ in $(M,$ $E),$ if all mixed cumulants of $M_{1}$ and $%
M_{2}$ vanish. The subsets $X_{1}$ and $X_{2}$ of $M$ are said to be free
over $N$ in $(M,$ $E),$ if the $W^{*}$-subalgebras $vN(X_{1},$ $N)$ and $%
vN(X_{2},$ $N)$ are free over $N$ in $(M,$ $E)$, where $vN(S_{1},$ $S_{2})$
is the von Neumann algebra generated by arbitrary sets $S_{1}$ and $S_{2}.$
In particular, we say that the $N$-valued random variables $x$ and $y$ are
free over $N$ in $(M,$ $E)$ if and only if $\{x\}$ and $\{y\}$ are free over 
$N$ in $(M,$ $E).$ Notice that the $N$-freeness is totally depending on the
conditional expectation $E.$ If $M_{1}$ and $M_{2}$ are free over $N$ in $%
(M, $ $E),$ then the $N$-free product von Neumann algebra $M_{1}$ $*_{N}$ $%
M_{2}$ is a $W^{*}$-subalgebra of $M,$ where

\strut

\begin{center}
$M_{1}*_{N}M_{2}=N\oplus \left( \oplus _{n=1}^{\infty }\left( \underset{%
i_{1}\neq i_{2},\,i_{2}\neq i_{3},\,...,\,i_{n-1}\neq i_{n}}{\oplus }%
(M_{i_{1}}^{o}\otimes ...\otimes M_{i_{n}}^{o})\right) \right) ,$
\end{center}

\strut \strut

where

\begin{center}
$M_{i_{j}}^{o}=M_{i_{j}}\ominus N,$ for all $j$ $=$ $1,$ ..., $n.$
\end{center}

\strut

Here, all algebraic operations $\oplus ,$ $\otimes $ and $\ominus $ are
defined under $W^{*}$-topology.

\strut

Also, if $(M_{1},$ $E_{1})$ and $(M_{2},$ $E_{2})$ are $N$-valued $W^{*}$%
-probability space with their conditional expectation $E_{j}$ $:$ $M_{j}$ $%
\rightarrow $ $N,$ for $j$ $=$ $1,$ $2.$ Then we can construct the free
product conditional expectation $E$ $=$ $E_{1}$ $*$ $E_{2}$ $:$ $M_{1}$ $%
*_{N}$ $M_{2}$ $\rightarrow $ $N$ making its cumulant $k_{n}^{(E)}(...)$
vanish for mixed $n$-tuples of $M_{1}$ and $M_{2}$ (See [19]). \strut

\strut \strut

The main result of this paper is that if $G_{1}$ $*$ $G_{2}$ is a free
product of groups $G_{1}$ and $G_{2},$ then

\strut

(0.5) $\ \ \ \ \ \ \ \ M$ $\times _{\alpha }(G_{1}$ $*$ $G_{2})$ $=$ $(M$ $%
\times _{\alpha }$ $G_{1})$ $*_{M}$ $(M$ $\times _{\alpha }$ $G_{2})$,

\strut

where $M$ is a von Neumann algebra and $\alpha $ $:$ $G_{1}$ $*$ $G_{2}$ $%
\rightarrow $ $AutM$ is an action. This shows that the group-freeness
implies a certain freeness on von Neumann algebras with amalgamation. Also,
this shows that, under the crossed product structure, the amalgamated
freeness determines the group freeness.

\strut

\begin{quote}
\textbf{Acknowledgment} \ The author really appreciates Prof. Florin
Radulescu for the valuable conversation and advice.
\end{quote}

\strut

\strut

\strut

\section{Crossed Product Probability Spaces}

\strut

\strut

In this chapter, we will observe some computations of $N$-valued moments and
cumulants of operators in the crossed product algebra $\Bbb{M}$ $=$ $M$ $%
\times _{\alpha }$ $G,$ with respect to the canonical conditional
expectation from $\Bbb{M}$ onto $M.$ Throughout this chapter, let $M$ be a
von Neumann algebra and $G,$ a group and let $\alpha $ $:$ $G$ $\rightarrow $
$AutM$ be an action of $G$ on $M,$ where $AutM$ is the automorphism group of 
$M.$

\strut

Denote the group identity of $G$ by $e_{G}.$ Consider the trivial subgroup $%
G_{0}$ $=$ $<e_{G}>$ of $G$ and the crossed product algebra $\Bbb{M}_{0}$ $=$
$M$ $\times _{\alpha }$ $G_{0}.$ Then this algebra $\Bbb{M}_{0}$ is a $W^{*}$%
-subalgebra of $\Bbb{M}$ and it satisfies that

\strut

(1.1) \ \ \ \ \ \ \ \ \ \ \ \ \ \ \ \ \ \ \ \ \ \ \ \ \ \ \ \ \ \ \ \ \ $%
\Bbb{M}_{0}=M,$

\strut

where the equality ``$=$'' means ``$*$-isomorphic''. Indeed, there exists a
linear map sending $m$ $\in $ $M$ to $m$ $u_{e_{G}}$ in $\Bbb{M}_{0}.$ This
is the $*$-isomorphism from $M$ onto $\Bbb{M}_{0},$ since

\strut

(1.2)$\ \ \ \ \ \ \ \ \ \ m_{1}m_{2}\longmapsto \left\{ 
\begin{array}{ll}
(m_{1}m_{2})u_{e_{G}} & =m_{1}\alpha _{e_{G}}(m_{2})u_{e_{G}} \\ 
& =m_{1}u_{e_{G}}m_{2}u_{e_{G}}u_{e_{G}} \\ 
& =(m_{1}u_{e_{G}})(m_{2}u_{e_{G}}),
\end{array}
\right. $

\strut

for all $m_{1},$ $m_{2}$ $\in $ $M.$ The first equality of the above formula
holds, because $\alpha _{e_{G}}$ is the identity automorphism on $M$
satisfying that $\alpha _{e_{G}}(m)$ $=$ $m,$ for all $m$ $\in $ $M.$ Also,
the third equality holds, because $u_{e_{G}}$ $u_{e_{G}}$ $=$ $u_{e_{G}^{2}}$
$=$ $u_{e_{G}}$ on $G_{0}$ (and also on $G$).

\strut

\begin{proposition}
Let $G_{0}$ $=$ $<e_{G}>$ be the trivial subgroup of $G$ and let $\Bbb{M}_{0}
$ $=$ $M$ $\times _{\alpha }$ $G_{0}$ be the crossed product algebra, where $%
\alpha $ is the given action of $G$ on $M.$ Then the von Neumann algebra $%
\Bbb{M}_{0}$ and $M$ are $*$-isomorphic. i.e., $\Bbb{M}_{0}$ $=$ $M.$ $%
\square $
\end{proposition}

\strut

From now, we will identify $M$ and $\Bbb{M}_{0},$ as $*$-isomorphic von
Neumann algebras.

\strut

\begin{definition}
Let $\Bbb{M}$ $=$ $M$ $\times _{\alpha }$ $G$ be the given crossed product
algebra. Define a canonical conditional expectation $E_{M}$ $:$ $\Bbb{M}$ $%
\rightarrow $ $M$ by

\strut 

(1.3) $\ \ \ \ \ \ E_{M}\left( \underset{g\in G}{\sum }m_{g}u_{g}\right)
=m_{e_{G}},$ for all $\underset{g\in G}{\sum }m_{g}u_{g}\in \Bbb{M}.$
\end{definition}

\strut

By (0.4), we have $u_{e_{G}}m=\alpha _{e_{G}}(m)\,u_{e_{G}}$ $=$ $m$ $%
u_{e_{G}}.$ So, indeed, the $\Bbb{C}$-linear map $E$ is a conditional
expectation; By the very definition, $E$ is continous and

\strut

(i) \ $E_{M}\left( m\right) =E_{M}\left( mu_{e_{G}}\right) =E_{M}\left(
u_{e_{G}}m\right) =m,$ for all $m$ $\in $ $M,$

\strut

(ii) $E_{M}\left( m_{1}(mu_{g})m_{2}\right) =m_{1}E_{M}(m_{g}u_{g})m_{2}$

$\strut $

\begin{center}
$=\left\{ 
\begin{array}{ll}
m_{1}m_{g}m_{2}=m_{1}E_{M}(mu_{g})m_{2} & \text{if }g=e_{G} \\ 
&  \\ 
0_{M}\text{ }=m_{1}E_{M}(mu_{g})m_{2} & \text{otherwise,}
\end{array}
\right. $
\end{center}

\strut

for all $m_{1},$ $m_{2}\in M$ and $mu_{g}\in \Bbb{M}.$ Therefore, we can
conclude that

\strut

\begin{center}
$E_{M}\left( m_{1}xm_{2}\right) =m_{1}E_{M}(x)m_{2},$ for $m_{1},$ $m_{2}$ $%
\in M$ and $x$ $\in $ $\Bbb{M}.$
\end{center}

\strut

(iii) For $\underset{g\in G}{\sum }m_{g}u_{g}\in \Bbb{M},$

$\strut $

\begin{center}
$
\begin{array}{ll}
E_{M}\left( (\underset{g\in G}{\sum }m_{g}u_{g})^{*}\right)  & 
\,=E_{M}\left( \underset{g\in G}{\sum }u_{g}^{*}m_{g}^{*}\right)  \\ 
&  \\ 
& 
\begin{array}{l}
=E_{M}\left( \underset{g\in G}{\sum }\alpha _{g}(m_{g}^{*})u_{g^{-1}}\right) 
\\ 
\\ 
=\alpha _{e_{G}}(m_{e_{G}}^{*})=m_{e_{G}}^{*}=E_{M}\left( \underset{g\in G}{%
\sum }m_{g}u_{g}\right) ^{*},
\end{array}
\end{array}
$
\end{center}

\strut

Therefore, by (i), (ii) and (iii), the map $E_{M}$ is a conditional
expectation. Thus the pair $(\Bbb{M},$ $E_{M})$ is a $M$-valued $W^{*}$%
-probability space.

\strut \strut \strut

\begin{definition}
The $M$-valued $W^{*}$-probability space $\left( \Bbb{M},\text{ }%
E_{M}\right) $ is called the $M$-valued crossed product probability space.
\end{definition}

\strut

It is trivial that $\Bbb{C}\cdot 1_{M}$ is a $W^{*}$-subalgebra of $M.$
Consider the crossed product $\Bbb{M}_{G}$ $=$ $\Bbb{C}$ $\times _{\alpha }$ 
$G,$ as a $W^{*}$-subalgebra of $\Bbb{M}.$ Recall the group von Neumann
algebra $L(G)$ defined by

\strut

\begin{center}
$L(G)=\overline{\Bbb{C}[G]}^{w}.$
\end{center}

\strut

Since every element $y$ in $\Bbb{M}_{G}$ has its Fourier expansion $y$ $=$ $%
\underset{g\in G}{\sum }$ $t_{g}u_{g}$ and since every element in $L(G)$ has
its Fourier expansion $\underset{g\in G}{\sum }$ $r_{g}u_{g},$ there exists
a $*$-isomorphism, which is the generator-preserving linear map, between $%
\Bbb{M}_{G}$ and $L(G).$

\strut

\begin{proposition}
Let $\Bbb{M}_{G}$ $\equiv $ $\Bbb{C}\cdot 1_{M}\times _{\alpha }G$ be the
crossed product algebra. Then $\Bbb{M}_{G}$ $=$ $L(G).$ $\square $
\end{proposition}

\strut \strut

\strut

\section{Moments and Cumulants on $(\Bbb{M},$ $E_{M})$}

\strut

\strut

In the previous section, we defined an amalgamated $W^{*}$-probability space
for the given crossed product algebra $\Bbb{M}$ $=$ $M$ $\times _{\alpha }$ $%
G.$ Throughout this chapter, we will let $M$ be a von Neumann algebra and $G,
$ a group and let $\alpha $ $:$ $G$ $\rightarrow $ $AutM$ be an action of $G$
on $M$. We will compute the amalgamated moments and cumulnats of operators
in $\Bbb{M}.$ These computations will play a key role to get our main
results (0.5), in Chapter 3.\strut  Let $(\Bbb{M},$ $E_{M})$ be the $M$%
-valued crossed product probability space.

\strut

\begin{quote}
\textbf{Notation} From now, we denote $\alpha _{g}(m)$ by $m^{g},$ for
convenience. $\square $
\end{quote}

\strut \strut

Consider group von Neumann algebras $L(G),$ which are $*$-isomorphic to $%
\Bbb{M}_{G}$ $=$ $\Bbb{C}$ $\times _{\alpha }$ $G,$ with its canonical trace 
$tr$ on it. On $L(G),$ we can always define its canonical trace $tr$ as
follows,

\strut

(2.1) $\ \ \ \ \ \ \ \ \ \ \ tr\left( \underset{g\in G}{\sum }%
r_{g}u_{g}\right) =r_{e_{G}},$ for all $\underset{g\in G}{\sum }%
r_{g}u_{g}\in L(G),$

\strut

where $r_{g}\in \Bbb{C}$, for $g$ $\in $ $G.$ So, the pair $(L(G),$ $tr)$ is
a $\Bbb{C}$-valued $W^{*}$-probability space. We can see that the unitary
representations $\{u_{g}\}_{g\in G}$ in $(\Bbb{M}$, $E)$ and $%
\{u_{g}\}_{g\in G}$ in $(L(G),$ $tr)$ are identically distributed.

\strut

By using the above new notation, we have

\strut

$\ \ \ \ \ \ \left( m_{g_{1}}u_{g_{1}}\right) \left(
m_{g_{2}}u_{g_{2}}\right) ...\left( m_{g_{n}}u_{g_{n}}\right) $

$\strut $

(2.2) $\ \ \ \ \ \ \ \ \ \ \ \ \ \ \ \ \ \ \ =\left(
m_{g_{1}}m_{g_{2}}^{g_{1}}m_{g_{3}}^{g_{1}g_{2}}...m_{g_{n}}^{g_{1}g_{2}...g_{n-1}}\right) u_{g_{1}...g_{n}}, 
$

\strut

for all $m_{g_{j}}u_{g_{j}}$ $\in $ $\Bbb{M},$ $j$ $=$ $1,$ ..., $n,$ where $%
n$ $\in $ $\Bbb{N}.$ The following lemma shows us that a certain collection
of $M$-valued random variables in $(\Bbb{M},$ $E_{M})$ and the generators of
group von Neumann algebra $(L(G),$ $tr)$ are identically distributed (over $%
\Bbb{C}$).

\strut \strut

\begin{lemma}
Let $u_{g_{1}},$ ..., $u_{g_{n}}$ $\in $ $\Bbb{M}$ (i.e., $u_{g_{k}}$ $=$ $%
1_{M}\cdot u_{g_{k}}$ in $\Bbb{M}$, for $k$ $=$ $1,$ ..., $n$.). Then

\strut 

(2.3) $\ \ \ \ \ \ \ \ \ \ \ \ \ \ \ E_{M}\left(
u_{g_{1}}...u_{g_{n}}\right) =tr\left( u_{g_{1}}...u_{g_{n}}\right) \cdot
1_{M},$

\strut 

where $tr$ is the canonical trace on the group von Neumann algebra $L(G).$
\end{lemma}

\strut

\begin{proof}
By definition of $E_{M},$

\strut

$\ E_{M}\left( u_{g_{1}}...u_{g_{m}}\right) =E_{M}\left( (1_{M}\cdot
1_{M}^{g_{1}}\cdot 1_{M}^{g_{1}g_{2}}\cdot \cdot \cdot
1_{M}^{g_{1}g_{2}...g_{n-1}})u_{g_{1}...g_{n}}\right) $

\strut

$\ \ \ \ \ \ \ \ \ \ \ \ \ \ \ \ \ \ \ \ \ =E_{M}\left(
u_{g_{1}...g_{n}}\right) $

\strut

since $%
1_{M}^{g}=u_{g}1_{M}u_{g^{-1}}=u_{g}u_{g^{-1}}=u_{gg^{-1}}=u_{e_{G}}=1_{M},$
for all $g$ $\in $ $G$

\strut

$\ \ \ \ \ \ \ \ \ \ \ \ \ \ \ \ \ \ \ \ \ =\left\{ 
\begin{array}{ll}
1_{M} & \text{if }g_{1}...g_{n}=e_{G} \\ 
0_{M} & \text{otherwise,}
\end{array}
\right. $

\strut

for all $n$ $\in $ $\Bbb{N}.$ By definition of $tr$ on $L(G),$ we have that

\strut

$\ \ \ \ \ \ tr\left( u_{g_{1}}...u_{g_{n}}\right) =tr\left(
u_{g_{1}...g_{n}}\right) =\left\{ 
\begin{array}{ll}
1 & \text{if }g_{1}...g_{n}=e_{G} \\ 
0 & \text{otherwise,}
\end{array}
\right. $

\strut

for all $n$ $\in $ $\Bbb{N}.$
\end{proof}

\strut

We want to compute the $M$-valued cumulant $k_{n}^{E_{M}}\left(
m_{g_{1}}u_{g_{1}},\text{ ..., }m_{g_{n}}u_{g_{n}}\right) ,$ for all $%
m_{g_{k}}u_{g_{k}}$ $\in $ $\Bbb{M}$ and $n$ $\in $ $\Bbb{N}.$ If this $M$%
-valued cumulant has a ``good'' relation with the cumulant $k_{n}^{tr}\left(
u_{g_{1}},\text{ ..., }u_{g_{n}}\right) ,$ then we might find the relation
between a group free product in $G$ and $M$-valued free product in $\Bbb{M}.$
The following three lemmas are the preparation for computing the $M$-valued
cumulant $k_{n}^{E_{M}}$ $(m_{g_{1}}u_{g_{1}},$ ..., $m_{g_{n}}u_{g_{n}}).$

\strut

\begin{lemma}
Let $(\Bbb{M},$ $E_{M})$ be the $M$-valued crossed product probability space
and let $m_{g_{1}}u_{g_{1}},$ ..., $m_{g_{n}}u_{g_{n}}$ be $M$-valued random
variables in $(\Bbb{M},$ $E_{M}),$ for $n$ $\in $ $\Bbb{N}.$ Then

\strut 

$\ \ \ \ \ \ E_{M}\left( m_{g_{1}}u_{g_{1}}...m_{g_{n}}u_{g_{n}}\right) $

\strut 

(2.4) $\ \ \ \ \ \ \ \ \ \ \ =\left\{ 
\begin{array}{ll}
m_{g_{1}}m_{g_{2}}^{g_{1}}m_{g_{3}}^{g_{1}g_{2}}...m_{g_{n}}^{g_{1}...g_{n-1}}
& \text{if }g_{1}...g_{n}=e_{G} \\ 
&  \\ 
0_{M} & \text{otherwise,}
\end{array}
\right. $

in $M.$
\end{lemma}

\strut

\begin{proof}
By the straightforward computation, we can get that

\strut

$\ \ \ \ \ \ E_{M}\left( m_{g_{1}}u_{g_{1}}...m_{g_{n}}u_{g_{n}}\right) $

\strut

$\ \ \ \ \ \ \ \ \ \ \ =E_{M}\left(
m_{g_{1}}m_{g_{2}}^{g_{1}}m_{g_{3}}^{g_{1}g_{2}}...m_{g_{n}}^{g_{1}...g_{n-1}}\cdot u_{g_{1}}u_{g_{2}}...u_{g_{n}}\right) 
$

\strut

by (0.2)

\strut

$\ \ \ \ \ \ \ \ \ \ \ =E_{M}\left(
(m_{g_{1}}m_{g_{2}}^{g_{1}}m_{g_{3}}^{g_{1}g_{2}}...m_{g_{n}}^{g_{1}...g_{n-1}})u_{g_{1}...g_{n}}\right) 
$

\strut \strut

$\ \ \ \ \ \ \ \ \ \ =\left(
m_{g_{1}}m_{g_{2}}^{g_{1}}m_{g_{3}}^{g_{1}g_{2}}...m_{g_{n}}^{g_{1}...g_{n-1}}\right) E_{M}\left( u_{g_{1}...g_{n}}\right) 
$

\strut

since $E_{M}$ $:$ $\Bbb{M}$ $\rightarrow $ $M$ $=$ $M$ $\times _{\alpha }$ $%
<e_{G}>$ is a conditional expectation

\strut

$\ \ \ \ \ \ \ \ \ \ =\left\{ 
\begin{array}{ll}
m_{g_{1}}m_{g_{2}}^{g_{1}}...m_{g_{n}}^{g_{1}...g_{n-1}} & \text{if }%
g_{1}...g_{n}=e_{G} \\ 
&  \\ 
0_{M} & \text{otherwise,}
\end{array}
\right. $

\strut

by the previous lemma. \strut \strut 
\end{proof}

\strut \strut 

Based on the previous lemma, we will compute the partition-depending moments
of $M$-valued random variables. But first, we need the following observation.

\strut

\begin{lemma}
Let $mu_{g}\in (\Bbb{M},$ $E_{M})$ be a $M$-valued random variable. Then $%
E_{M}\left( u_{g}m\right) =m^{g}E_{M}\left( u_{g}\right) .$
\end{lemma}

\strut

\begin{proof}
Compute

\strut

$\ \ \ \ \ \ \ E_{M}\left( u_{g}m\right) =E_{M}\left(
u_{g}mu_{g^{-1}}u_{g}\right) =E_{M}\left( m^{g}u_{g}\right)
=m^{g}E_{M}\left( u_{g}\right) .$

\strut
\end{proof}

\strut

Since $E_{M}$ is a conditional expectation, $E_{M}\left( u_{g}m\right) $ $=$ 
$E_{M}\left( u_{g}\right) $ $m,$ too. So, by the previous lemma, we have that

\strut

(2.6) $\ \ \ \ \ \ \ \ \ E_{M}\left( u_{g}\right) m=E\left( u_{g}m\right)
=m^{g}E\left( u_{g}\right) .$

\strut \strut \strut

In \ the following lemma, we will extend this observation (2.6) to the
general case. Notice that since $E_{M}$ is a $M$-valued conditional
expectation, we have to consider the insertion property (See [19]). i.e., in
general,

\strut

\begin{center}
$E_{M,\pi }\left( x_{1},...,x_{n}\right) \neq $ $\underset{V\in \pi }{\Pi }%
E_{M,V}\left( x_{1},...,x_{n}\right) ,$
\end{center}

\strut

for $x_{1},$ ..., $x_{n}$ $\in $ $\Bbb{M}$, where $E_{M,V}(...)$ is the
block-depending moments. But, if $x_{k}$ $=$ $u_{g_{k}}$ $=$ $1_{M}\cdot
u_{g_{k}}$ in $\Bbb{M},$ then we can have that

\strut

\begin{center}
$E_{M,\text{ }\pi }\left( u_{g_{1}},...,u_{g_{n}}\right) =\underset{B\in \pi 
}{\Pi }E_{M,V}\left( u_{g_{1}},...,u_{g_{n}}\right) $
\end{center}

\strut

since

\begin{center}
$E_{M}\left( u_{g}\right) $ $=$ $\left\{ 
\begin{array}{lll}
1\in \Bbb{C}\text{ }\cdot 1_{M} &  & \text{if }g=e_{G} \\ 
0\in \Bbb{C}\text{ }\cdot 1_{M} &  & \text{otherwise,}
\end{array}
\right. $
\end{center}

and hence

\strut

\begin{center}
$=\underset{B\in \pi }{\Pi }\left( tr_{V}\left(
u_{g_{1}},...,u_{g_{n}}\right) \cdot 1_{M}\right) $
\end{center}

\strut

by (2.3).

\strut

Suppose that $\pi $ $\in $ $NC(n)$ is a partition which is not $1_{n}$ and
by [$V$ $\in $ $\pi $], denote the relation [$V$ is a block of $\pi $]. We
say that a block $V$ $=$ $(j_{1},$ ..., $j_{p})$ is inner in a block $B$ $=$ 
$(i_{1},$ ..., $i_{k})$, where $V,$ $B$ $\in $ $\pi $, if there exists $k_{0}
$ $\in $ $\{2,$ ..., $k$ $-$ $1\}$ such that $i_{k_{0}}$ $<$ $j_{t}$ $<$ $%
i_{k_{0}+1},$ for all $t$ $=$ $1,$ ..., $p.$ In this case, we also say that $%
B$ is outer than $V.$ Also, we say that $V$ is innerest if there is no other
block inner in $V.$ For instance, if we have a partition

\strut

\begin{center}
$\pi $ $=$ $\{(1,6),(2,5),(3,4)\}$ in $NC(6).$
\end{center}

\strut \strut

Then the block $(2,$ $5)$ is inner in the block $(1,$ $6)$ and the block $(3,
$ $4)$ is inner in the block $(2,$ $5).$ Clearly, the block $(3,$ $4)$ is
inner in both $(2,5)$ and $(1,$ $6),$ and there is no other block inner in $%
(3,$ $4).$ So, the block $(3,$ $4)$ is an innerest block in $\pi .$ Remark
that it is possible there are several innerest blocks in a certain
noncrossing partition. Also, notice that if $V$ is an innerest block, then
there exists $j$ such that $V$ $=$ $(j,$ $j$ $+$ $1,$ ..., $j$ $+$ $\left|
V\right| $ $-$ $1),$ where $\left| V\right| $ means the cardinality of
entries of $V.$

\strut

\begin{lemma}
Let $n$ $\in $ $\Bbb{N}$ and $\pi $ $\in $ $NC(n),$ and let $%
m_{g_{1}}u_{g_{1}},$ ..., $m_{g_{n}}u_{g_{n}}$ $\in $ $(\Bbb{M}$, $E_{M})$
be the $M$-valued random variables. Then

\strut \strut \strut 

$\ \ \ \ \ \ E_{M,\pi }\left( m_{g_{1}}u_{g_{1}},\text{ ..., }%
m_{g_{n}}u_{g_{n}}\right) $

\strut 

(2.7) $\ \ \ \ \ \ \ \ \ \ \ =\left(
m_{g_{1}}m_{g_{2}}^{g_{1}}...m_{g_{n}}^{g_{1}...g_{n-1}}\right) tr_{\pi
}\left( u_{g_{1}},...,u_{g_{n}}\right) ,$

\strut 

where $tr$ is the canonical trace on the group von Neumann algebra $L(G).$
\end{lemma}

\strut

\begin{proof}
If $\pi $ $=$ $1_{n},$ then $E_{M,1_{n}}(...)$ $=$ $E_{M}(...)$ and $%
tr_{1_{n}}(...)$ $=$ $tr(...),$ and hence we are done, by (2.3) and (2.4).
Assume that $\pi $ $\neq $ $1_{n}$ in $NC(n).$ Assume that $V$ $=$ $(j,$ $j$ 
$+$ $1,$ ..., $j$ $+$ $k)$ is an innerest block of $\pi $. Then

\strut

$\ \ \ \ \ T_{V}\overset{def}{=}E_{M,\text{ }V}\left(
m_{g_{1}}u_{g_{1}},...,m_{g_{n}}u_{g_{n}}\right) $

\strut

$\ \ \ \ \ \ \ \ \ \ =E_{M}\left(
m_{g_{j}}u_{g_{j}}m_{g_{j+1}}u_{g_{j+1}}...m_{g_{j+k}}u_{g_{j+k}}\right) $

\strut

$\ \ \ \ \ \ \ \ \ \ =\left(
m_{g_{j}}m_{g_{j+1}}^{g_{j}}m_{g_{j+2}}^{g_{j}g_{j+1}}...m_{g_{j+k}}^{g_{j}g_{j+1}...g_{j+k-1}}\right) \cdot tr\left( u_{g_{j}...g_{j+k}}\right) . 
$

\strut

Suppose $V$ is inner in a block $B$ of $\pi $ and $B$ is inner in all other
blocks $B^{\prime }$ where $V$ is inner in $B^{\prime }.$ Let $B$ $=$ $%
(i_{1},$ ..., $i_{k})$ and assume that there is $k_{0}$ $\in $ $\{2,$ ..., $%
k $ $-$ $1\}$ such that $i_{k_{0}}$ $<$ $t$ $<$ $i_{k_{0}+1},$ for all $t$ $%
= $ $j,$ $j$ $+$ $1,$ ..., $j$ $+$ $k.$ Then the $B$-depending moment goes to

\strut

$\ E_{M}\left(
m_{g_{i_{1}}}u_{g_{i_{1}}}...m_{g_{k_{0}}}u_{g_{k_{0}}}(T_{V})m_{g_{k_{0}+1}}u_{g_{k_{0}+1}}...m_{g_{i_{k}}}u_{g_{i_{k}}}\right) 
$

\strut

$\ =E_{M}\left( \left(
m_{g_{i_{1}}}m_{g_{i_{2}}}^{g_{i_{1}}}...m_{g_{k_{0}}}^{g_{i_{1}}...g_{i_{k_{0}-1}}}\right. \right. 
$

\strut

$\ \ \ \ \ \ \ \ \ \ \ \cdot \left(
m_{g_{j}}^{g_{i_{1}}...g_{i_{k_{0}}}}m_{g_{j+1}}^{g_{i_{1}}...g_{i_{k_{0}}}g_{j}}...m_{g_{j+1}}^{g_{i_{1}}...g_{i_{k_{0}}}g_{j}..g_{j+k-1}}\right) 
$

\strut

$\ \ \ \ \ \ \ \ \ \left. \left. \cdot
m_{g_{i_{k_{0}}}}^{g_{i_{1}}...g_{i_{k_{0}}}g_{j}...g_{j+k}}...m_{g_{i_{k}}}^{g_{i_{1}}...g_{j}...g_{j+1}...g_{i_{k-1}}}\right) u_{g_{i_{1}}...g_{i_{k_{0}}}g_{j}...g_{j+k}g_{i_{k_{0}+1}}...g_{i_{k}}}\right) 
$

\strut

$\ =\left(
m_{g_{i_{1}}}m_{g_{i_{2}}}^{g_{i_{1}}}...m_{g_{k_{0}}}^{g_{i_{1}}...g_{i_{k_{0}-1}}}\right. 
$

\strut

$\ \ \ \ \ \ \ \ \cdot \left(
m_{g_{j}}^{g_{i_{1}}...g_{i_{k_{0}}}}m_{g_{j+1}}^{g_{i_{1}}...g_{i_{k_{0}}}g_{j}}...m_{g_{j+1}}^{g_{i_{1}}...g_{i_{k_{0}}}g_{j}..g_{j+k-1}}\right) 
$

\strut

$\ \ \ \ \ \ \ \ \left. \cdot
m_{g_{i_{k_{0}}}}^{g_{i_{1}}...g_{i_{k_{0}}}g_{j}...g_{j+k}}...m_{g_{i_{k}}}^{g_{i_{1}}...g_{j}...g_{j+1}...g_{i_{k-1}}}\right) E_{M}\left( u_{g_{i_{1}}...g_{i_{k_{0}}}g_{j}...g_{j+k}g_{i_{k_{0}+1}}...g_{i_{k}}}\right) . 
$

\strut

By doing the above process for all block-depending moments in the $\pi $%
-depending moments, we can get that

\strut

$\ \ \ E_{M,\pi }\left( m_{g_{1}}u_{g_{1}},\text{ ..., }m_{g_{n}}u_{g_{n}}%
\right) $

\strut

$\ \ \ \ \ \ \ \ \ \ \ \ \ \ =\left(
m_{g_{1}}m_{g_{2}}^{g_{1}}m_{g_{3}}^{g_{1}g_{2}}...m_{g_{n}}^{g_{1}...g_{n-1}}\right) E_{\pi }\left( u_{g_{1}},%
\text{ ..., }u_{g_{n}}\right) .$

\strut

By (2.3), we know $E_{\pi }\left( u_{g_{1}},...,u_{g_{n}}\right) =$ $tr_{\pi
}(u_{g_{1}},$ ..., $u_{g_{n}})$ $\cdot $ $1_{M},$ where $tr$ is the
canonical trace on the group von Neumann algebra $L(G).$
\end{proof}

\strut

By the previous lemmas and proposition, we have the following theorem.

\strut

\begin{theorem}
Let $m_{g_{1}}u_{g_{1}},$ ..., $m_{g_{n}}u_{g_{n}}$ $\in $ $(\Bbb{M},$ $%
E_{M})$ be the $M$-valued random variables, for $n$ $\in $ $\Bbb{N}.$ Then

\strut 

$\ \ \ k_{n}^{E_{M}}\left( m_{g_{1}}u_{g_{1}},\text{ ..., }%
m_{g_{n}}u_{g_{n}}\right) $

\strut 

(2.8) $\ \ \ \ \ \ =\left(
m_{g_{1}}m_{g_{2}}^{g_{1}}m_{g_{3}}^{g_{1}g_{2}}...m_{g_{n}}^{g_{1}...g_{n-1}}\right) k_{n}^{tr}\left( u_{g_{1}},...,u_{g_{n}}\right) .
$
\end{theorem}

\strut

\begin{proof}
Observe that

\strut

$\ \ \ k_{n}^{M}\left( m_{g_{1}}u_{g_{1}},...,m_{g_{n}}u_{g_{n}}\right) $

\strut

$\ \ \ =\underset{\pi \in NC(n)}{\sum }\,E_{M,\text{ }\pi }\left(
m_{g_{1}}u_{g_{1}},...,m_{g_{n}}u_{g_{n}}\right) \mu (\pi ,1_{n})$

\strut

$\ \ \ =\underset{\pi \in NC(n)}{\sum }\left(
(m_{g_{1}}m_{g_{2}}^{g_{1}}...m_{g_{n}}^{g_{1}...g_{n-1}})\,tr_{\pi }\left(
u_{g_{1}},...,u_{g_{n}}\right) \right) \mu (\pi ,1_{n})$

\strut

by (2.7)

\strut

$\ \ \ =\left(
m_{g_{1}}m_{g_{2}}^{g_{1}}m_{g_{3}}^{g_{1}g_{2}}...m_{g_{n}}^{g_{1}...g_{n-1}}\right) \left( 
\underset{\pi \in NC(n)}{\sum }tr_{\pi }(u_{g_{1}},...,u_{g_{n}})\mu (\pi
,1_{n})\right) $

\strut

$\ \ \ =\left(
m_{g_{1}}m_{g_{2}}^{g_{1}}m_{g_{3}}^{g_{1}g_{2}}...m_{g_{n}}^{g_{1}...g_{n-1}}\right) \,k_{n}^{tr}\left( u_{g_{1}},...,u_{g_{n}}\right) . 
$

\strut
\end{proof}

\strut

The above theorem shows us that there is close relation between the $M$%
-valued cumulant on $(\Bbb{M},$ $E_{M})$ and $\Bbb{C}$-valued cumulant on $%
(L(G),$ $tr).$

\strut

\begin{example}
In this example, instead of using (2.7) directly, we will compute the $\pi $%
-depending moment of $m_{g_{1}}$ $u_{g_{1}},$ ..., $m_{g_{n}}$ $u_{g_{n}}$
in $\Bbb{M},$ only by using the simple computations (0.1) \symbol{126}
(0.4). By doing this, we can understand why (2.7) holds concretely. Let $\pi 
$ $=$ $\{(1,$ $4),$ $(2,$ $3),$ $(5)\}$ in $NC(5).$ Then

\strut 

$\ E_{M,\,\pi }\left(
m_{g_{1}}u_{g_{1}},\,\,...\,\,,m_{g_{5}}u_{g_{5}}\right) $

\strut 

$\ \ =E_{M}\left( m_{g_{1}}u_{g_{1}}E_{M}\left(
m_{g_{2}}u_{g_{2}}m_{g_{3}}u_{g_{3}}\right) m_{g_{4}}u_{g_{4}}\right)
E_{M}\left( m_{g_{5}}u_{g_{5}}\right) $

\strut 

$\ \ =m_{g_{1}}E_{M}\left( u_{g_{1}}E_{M}\left(
m_{g_{2}}m_{g_{3}}^{g_{2}}u_{g_{2}g_{3}}\right) m_{g_{4}}u_{g_{4}}\right)
\,\left( \,m_{g_{5}}E_{M}(u_{g_{5}})\right) $

\strut 

$\ \ =m_{g_{1}}E_{M}\left(
u_{g_{1}}(m_{g_{2}}m_{g_{3}}^{g_{2}})E_{M}(u_{g_{2}g_{3}})m_{g_{4}}u_{g_{4}}%
\right) \left( m_{g_{5}}E_{M}(u_{g_{5}})\right) $

\strut 

$\ \ =m_{g_{1}}E_{M}\left(
m_{g_{2}}^{g_{1}}m_{g_{3}}^{g_{1}g_{2}}u_{g_{1}}E_{M}(u_{g_{2}g_{3}})m_{g_{4}}u_{g_{4}}\right) \left( m_{g_{5}}E_{M}(u_{g_{5}})\right) 
$

\strut 

$\ \ =m_{g_{1}}E_{M}\left(
m_{g_{2}}^{g_{1}}m_{g_{3}}^{g_{1}g_{2}}u_{g_{1}}m_{g_{4}}^{g_{2}g_{3}}E_{M}(u_{g_{2}g_{3}})u_{g_{4}}\right) \left( m_{g_{5}}E_{M}(u_{g_{5}})\right) 
$

\strut 

$\ \ =m_{g_{1}}m_{g_{2}}^{g_{1}}m_{g_{3}}^{g_{1}g_{2}}E_{M}\left(
u_{g_{1}}m_{g_{4}}^{g_{2}g_{3}}E_{M}(u_{g_{2}g_{3}})u_{g_{4}}\right) \left(
m_{g_{5}}E_{M}(u_{g_{5}})\right) $

\strut 

$\ \ =m_{g_{1}}m_{g_{2}}^{g_{1}}m_{g_{3}}^{g_{1}g_{2}}E_{M}\left(
m_{g_{4}}^{g_{1}g_{2}g_{3}}u_{g_{1}}E_{M}(u_{g_{2}g_{3}})u_{g_{4}}\right)
\left( m_{g_{5}}E_{M}(u_{g_{5}})\right) $

\strut 

$\ \
=m_{g_{1}}m_{g_{2}}^{g_{1}}m_{g_{3}}^{g_{1}g_{2}}m_{g_{4}}^{g_{1}g_{2}g_{3}}E_{M}\left( u_{g_{1}}E_{M}(u_{g_{2}g_{3}})u_{g_{4}}\right) m_{g_{5}}\left( E_{M}(u_{g_{5}})\right) 
$

\strut 

$\ \
=m_{g_{1}}m_{g_{2}}^{g_{1}}m_{g_{3}}^{g_{1}g_{2}}m_{g_{4}}^{g_{1}g_{2}g_{3}}E_{M}\left( u_{g_{1}}E_{M}(u_{g_{2}g_{3}})u_{g_{4}}m_{g_{5}}\right) \left( E_{M}(u_{g_{5}})\right) 
$

\strut 

$\ \
=m_{g_{1}}m_{g_{2}}^{g_{1}}m_{g_{3}}^{g_{1}g_{2}}m_{g_{4}}^{g_{1}g_{2}g_{3}}E_{M}\left( u_{g_{1}}E_{M}(u_{g_{2}g_{3}})m_{g_{5}}^{g_{4}}u_{g_{4}}\right) \left( E_{M}(u_{g_{5}})\right) 
$

\strut 

$\ \
=m_{g_{1}}m_{g_{2}}^{g_{1}}m_{g_{3}}^{g_{1}g_{2}}m_{g_{4}}^{g_{1}g_{2}g_{3}}E_{M}\left( u_{g_{1}}m_{g_{5}}^{g_{2}g_{3}g_{4}}E_{M}(u_{g_{2}g_{3}})u_{g_{4}}\right) \left( E_{M}(u_{g_{5}})\right) 
$

\strut 

$\ \
=m_{g_{1}}m_{g_{2}}^{g_{1}}m_{g_{3}}^{g_{1}g_{2}}m_{g_{4}}^{g_{1}g_{2}g_{3}}E_{M}\left( m_{g_{5}}^{g_{1}g_{2}g_{3}g_{4}}u_{g_{1}}E_{M}(u_{g_{2}g_{3}})u_{g_{4}}\right) \left( E_{M}(u_{g_{5}})\right) 
$

\strut 

$\ \ =\left(
m_{g_{1}}m_{g_{2}}^{g_{1}}m_{g_{3}}^{g_{1}g_{2}}m_{g_{4}}^{g_{1}g_{2}g_{3}}m_{g_{5}}^{g_{1}g_{2}g_{3}g_{4}}\right) \left( \left( E_{M}(u_{g_{1}}E_{M}(u_{g_{2}}u_{g_{3}})u_{g_{4}})\right) \left( E_{M}(u_{g_{5}})\right) \right) 
$

\strut 

$\ \ =\left(
m_{g_{1}}m_{g_{2}}^{g_{1}}m_{g_{3}}^{g_{1}g_{2}}m_{g_{4}}^{g_{1}g_{2}g_{3}}m_{g_{5}}^{g_{1}g_{2}g_{3}g_{4}}\right) \left( tr\left( u_{g_{1}}tr(u_{g_{2}}u_{g_{3}})u_{g_{4}}\right) \left( tr(u_{g_{5}})\right) \right) 
$

\strut 

$\ \ =\left(
m_{g_{1}}m_{g_{2}}^{g_{1}}m_{g_{3}}^{g_{1}g_{2}}m_{g_{4}}^{g_{1}g_{2}g_{3}}m_{g_{5}}^{g_{1}g_{2}g_{3}g_{4}}\right) \,\left( tr_{\pi }(u_{g_{1}},u_{g_{2}},u_{g_{3}},u_{g_{4}},u_{g_{5}})\right) .
$

\strut 
\end{example}

\strut

\begin{example}
We can compute the following $M$-valued cumulant, by applying (2.8).

\strut 

$\ k_{3}^{E_{M}}\left(
m_{g_{1}}u_{g_{1}},m_{g_{2}}u_{g_{2}},m_{g_{3}}u_{g_{3}}\right) =\left(
m_{g_{1}}m_{g_{2}}^{g_{1}}m_{g_{3}}^{g_{1}g_{2}}\right) \cdot
k_{3}^{tr}\left( u_{g_{1}},u_{g_{2}},u_{g_{3}}\right) $

\strut 

$\ \ \ =\left( m_{g_{1}}m_{g_{2}}^{g_{1}}m_{g_{3}}^{g_{1}g_{2}}\right)
\left( tr\left( u_{g_{1}g_{2}g_{3}}\right)
-tr(u_{g_{1}})tr(u_{g_{2}}u_{g_{3}})\right. $

\strut 

$\ \ \ \ \ \ \ \ \ \ \ \ \ \ \ \ \ \ \ \ \ \ \ \ \ \ \ \left.
-tr(u_{g_{1}}u_{g_{2}})tr(u_{g_{3}})+2tr(u_{g_{1}})tr(u_{g_{2}})tr(u_{g_{3}})\right) .
$
\end{example}

\strut \strut \strut

\strut \strut \strut \strut \strut \strut

\strut

\section{The Main Result (0.5)}

\strut

\strut

In this chapter, we will prove our main result (0.5). Like before,
throughout this chapter, let $M$ be a von Neumann algebra and $G,$ a group
and let $\alpha $ $:$ $M$ $\rightarrow $ $AutM$ be an action of $G$ on $M$.
Assume that a group $G$ is a group free product $G_{1}$ $*$ $G_{2}$ of
groups $G_{1}$ and $G_{2}.$ (Also, we can assume that there is a subgroup $%
G_{1}$ $*$ $G_{2}$ in the group $G,$ and $M$ $\times _{\alpha }$ $(G_{1}$ $*$
$G_{2})$ is a $W^{*}$-subalgebra of $\Bbb{M}.$) Recall that, by Voiculescu,
it is well-known that

\strut

\begin{center}
$L\left( G_{1}*G_{2}\right) =L(G_{1})*L(G_{2}),$
\end{center}

\strut

where ``$*$'' in the left-hand side is the group free product and ``$*$'' in
the right-hand side is the von Neumann algebra free product, where $L(K)$ is
a group von Neumann algebra of an arbitrary group $K.$ This says that the $%
\Bbb{C}$-freeness on $(L(G),$ $tr)$ is depending on the group freeness on $G$
$=$ $G_{1}$ $*$ $G_{2}.$ In other words, if the groups $G_{1}$ and $G_{2}$
are free in $G$ $=$ $G_{1}$ $*$ $G_{2},$ then the group von Neumann algebras 
$L(G_{1})$ and $L(G_{2})$ are free in $(L(G),$ $tr)$. Also, if two group von
Neumann algebras $L(G_{1})$ and $L(G_{2})$ are given and if we construct the 
$\Bbb{C}$-free product $L(G_{1})$ $*$ $L(G_{2})$ of them, with respect to
the canonical trace $tr_{G}$ $=$ $tr_{G_{1}}$ $*$ $tr_{G_{2}},$ where $%
tr_{G_{k}}$ is the canonical trace on $L(G_{k}),$ for $k$ $=$ $1,$ $2,$ then
this $\Bbb{C}$-free product is $*$-isomorphic to a group von Neumann algebra 
$L(G),$ where $G$ is the group free product $G_{1}$ $*$ $G_{2}$ of $G_{1}$
and $G_{2}.$

\strut

\begin{theorem}
Let $\Bbb{M}$ $=$ $M$ $\times _{\alpha }$ $G$ be a crossed product algebra,
where $G$ $=$ $G_{1}$ $*$ $G_{2}$ is the group free product of $G_{1}$ and $%
G_{2}.$ Then

\strut 

(3.1) $\ \ \ \ \ \ \ \ \ \ \ \ \ \Bbb{M}=\left( M\times _{\alpha
}G_{1}\right) *_{M}\left( M\times _{\alpha }G_{2}\right) ,$

\strut \strut 

where ``$*_{M}$'' is the $M$-valued free product of von Neumann algebras.
\end{theorem}

\strut

\begin{proof}
Let $G$ $=$ $G_{1}$ $*$ $G_{2}$ be the group free product of $G_{1}$ and $%
G_{2}.$ By Chapter 1, the crossed product algebra $\Bbb{M}$ has its $W^{*}$%
-subalgebra $M$ $=$ $M\ \times _{\alpha }$ $<e_{G}>,$ where $<e_{G}>$ is the
trivial subgroup of $G$ generated by the group identity $e_{G}$ $\in $ $G.$
Define the canonical conditional expectation $E_{M}$ $:$ $\Bbb{M}$ $%
\rightarrow $ $M$ by

\strut

$\ \ \ \ \ \ \ \ E_{M}\left( \underset{g\in G}{\sum }m_{g}u_{g}\right)
=m_{e_{G}},$ \ for all $\ \underset{g\in G}{\sum }m_{g}u_{g}\in \Bbb{M}.$

\strut

By (2.8), if $m_{g_{1}}u_{g_{1}},$ ..., $m_{g_{n}}u_{g_{n}}$ $\in $ $(\Bbb{M}%
,$ $E_{M})$ are $M$-valued random variables, then

\strut

\ \ \ \ \ $k_{n}^{E_{M}}\left(
m_{g_{1}}u_{g_{1}},...,m_{g_{n}}u_{g_{n}}\right) $

\strut

$\ \ \ \ \ \ \ \ \ \ \ \ \ \ \ \ \ =\left(
m_{g_{1}}m_{g_{2}}^{g_{1}}m_{g_{3}}^{g_{1}g_{2}}...m_{g_{n}}^{g_{1}...g_{n-1}}\right) \,k_{n}^{tr}\left( u_{g_{1}},%
\text{ ..., }u_{g_{n}}\right) ,$

\strut

for all $n$ $\in $ $\Bbb{N},$ where $tr$ is the canonical trace on $L(G).$
As we mentioned in the previous paragraph, the $\Bbb{C}$-freeness on $L(G)$
is completely determined by the group freeness of $G_{1}$ and $G_{2}$ on $G$
and vice versa. By the previous cumulant relation, the $M$-freeness on $\Bbb{%
M}$ is totally determined by the $\Bbb{C}$-freeness on $L(G).$ Therefore,
the $M$-freeness on $\Bbb{M}$ is determined by the group freeness on $G.$
Thus, we can conclude that

\strut

$\ \ \ \ \ \ \ \ \ \ \ \ \ \ M\times _{\alpha }(G_{1}*G_{2})=\left( M\times
_{\alpha }G_{1}\right) *_{M}(M\times _{\alpha }G_{2}).$

\strut \strut \strut 
\end{proof}

\strut

If $F_{N}$ is the free group with $N$-generators, then $L(F_{N})$ $=$ $%
*_{k=1}^{N}L(\Bbb{Z})_{k},$ where $L(\Bbb{Z})_{k}$ $=$ $L(\Bbb{Z}),$ for all 
$k$ $=$ $1,$ ..., $N.$ Also, $L(F_{N})$ $=$ $L(F_{k_{1}})$ $*$ $%
L(F_{k_{2}}), $ for all $k_{1},$ $k_{2}$ $\in $ $\Bbb{C}$ such that $k_{1}$ $%
+$ $k_{2}$ $=$ $N.$

\strut

\begin{corollary}
Let $F_{N}$ be the free group with $N$-generators, for $N$ $\in $ $\Bbb{N}$.
Then

\strut 

(3.2) $\ \ \ \ \ \ \ \ \ M$ $\times _{\alpha }$ $F_{N}=$ $\ \underset{N\text{%
-times}}{\underbrace{\left( M\times _{\alpha }\Bbb{Z}\right)
*_{M}...*_{M}(M\times _{\alpha }\Bbb{Z})}}$

and

\strut 

(3.3) \ \ $\ \ \ \ \ \ \ \ M\times _{\alpha }F_{N}=\left( M\times _{\alpha
}F_{k_{1}}\right) *_{M}\left( M\times _{\alpha }F_{k_{2}}\right) ,$

\strut 

whenever $k_{1}+k_{2}=N,$ for $k_{1},$ $k_{2}$ $\in $ $\Bbb{N}.$ $\ \square $
\end{corollary}

\strut \strut \strut \strut \strut

\strut \strut \strut \strut

\strut \strut

\strut \textbf{References}

\strut

\strut

\begin{quote}
{\small [1] \ \ A. G. Myasnikov and V. Shapilrain (editors), Group Theory,
Statistics and Cryptography, Contemporary Math, 360, (2003) AMS.}

{\small [2] \ \ A. Nica, R-transform in Free Probability, IHP course note,
available at www.math.uwaterloo.ca/\symbol{126}anica.}

{\small [3]\strut \ \ \ A. Nica and R. Speicher, R-diagonal Pair-A Common
Approach to Haar Unitaries and Circular Elements, (1995), www
.mast.queensu.ca/\symbol{126}speicher.\strut }

{\small [4] \ \ A. Nica, D. Shlyakhtenko and R. Speicher, R-cyclic Families
of Matrices in Free Probability, J. of Funct Anal, 188 (2002), 227-271.}

{\small [5] \ }$\ ${\small B. Solel, You can see the arrows in a Quiver
Operator Algebras, (2000), preprint.\strut }

{\small [6] \ \ D. Shlyakhtenko, Some Applications of Freeness with
Amalgamation, J. Reine Angew. Math, 500 (1998), 191-212.\strut }

{\small [7] \ \ D.Voiculescu, K. Dykemma and A. Nica, Free Random Variables,
CRM Monograph Series Vol 1 (1992).\strut }

{\small [8] \ \ D. Voiculescu, Operations on Certain Non-commuting
Operator-Valued Random Variables, Ast\'{e}risque, 232 (1995), 243-275.\strut 
}

{\small [9] \ \ D. Shlyakhtenko, A-Valued Semicircular Systems, J. of Funct
Anal, 166 (1999), 1-47.\strut }

{\small [10]\ D.W. Kribs and M.T. Jury, Ideal Structure in Free Semigroupoid
Algebras from Directed Graphs, preprint}

{\small [11]\ D.W. Kribs and S.C. Power, Free Semigroupoid Algebras, preprint%
}

{\small [12]\ G. C. Bell, Growth of the Asymptotic Dimension Function for
Groups, (2005) Preprint.}

{\small [13]\ I. Cho, Random Variables in a Graph }$W^{*}${\small %
-Probability Space, (2005) Ph. D. Thesis, Univ. of Iowa.}

{\small [14]\ I. Cho, Moments of Block Operators of a Group von Neumann
Algebra, (2005) Submitted to Manu.Math.}

{\small [15]\ I. Cho, Moments of the Radical Operator of a Group von Neumann
Algebra, (2005), Preprint. }

{\small [16] J. Stallings, Centerless Groups-An Algebraic Formulation of
Gottlieb's Theorem, Topology, Vol 4, (1965) 129 - 134.}

{\small [17] P.\'{S}niady and R.Speicher, Continous Family of Invariant
Subspaces for R-diagonal Operators, Invent Math, 146, (2001) 329-363.}

{\small [18] R. Gliman, V. Shpilrain and A. G. Myasnikov (editors),
Computational and Statistical Group Theory, Contemporary Math, 298, (2001)
AMS.}

{\small [19] R. Speicher, Combinatorial Theory of the Free Product with
Amalgamation and Operator-Valued Free Probability Theory, AMS Mem, Vol 132 ,
Num 627 , (1998).}

{\small [20] R. Speicher, Combinatorics of Free Probability Theory IHP
course note, available at www.mast.queensu.ca/\symbol{126}speicher.}

{\small [21] V. Jones, Subfactor and Knots }

{\small [22] F. Radulescu, Random Matrices, Amalgamated Free Products and
Subfactors of the von Neumann Algebra of a Free Group, of Noninteger Index,
Invent. Math., 115, (1994) 347 - 389.}
\end{quote}

\end{document}